\documentclass[12pt]{article}

\newcommand{\bc}{\begin{center}}
\newcommand{\ec}{\end{center}}
\newcommand{\be}{\begin{equation}}
\newcommand{\ee}{\end{equation}}
\newcommand{\ba}{\begin{array}}
\newcommand{\ea}{\end{array}}
\newcommand{\bea}{\begin{eqnarray}}
\newcommand{\eea}{\end{eqnarray}}
\newcommand{\edc}{\end{document}}

\def\l{\lambda}

\begin{document}
\thispagestyle{empty}
\begin{center}

{\bf ON $F$-QUADRATIC STOCHASTIC OPERATORS}\\
\vspace{0.5cm}

{\bf U.A. Rozikov$^1$, U.U. Jamilov$^2$}\\[0.4mm]
{$^1$\it Institute of Mathematics, Tashkent, Uzbekistan,} e-mail: rozikovu@yandex.ru\\
{$^2$\it College n.43, Romitan, Bukhara region, Uzbekistan.} e-mail: jamilovu@yandex.ru\\
\end{center}

 {\bf Abstract.}
 In this paper we introduce a notion of $F-$ quadratic stochastic
operator. For a wide class of such operators we show that each
operator of the class has unique fixed point. Also we prove that any
trajectory of the $F$-quadratic stochastic operator  converges to
the fixed point exponentially fast.


\section{Introduction}

 Quadratic stochastic operator frequently arise in many models of mathematical
genetics [1, 6-8], [14],[15].

The quadratic stochastic operator (QSO) is a mapping of the simplex
$$ S^{m-1}=\{x=(x_1,...,x_m)\in R^m: x_i\geq 0, \sum^m_{i=1}x_i=1 \} \eqno(1)$$
into itself, of the form
$$ V: x_k'=\sum^n_{i,j=1}p_{ij,k}x_ix_j, \ \ (k=1,...,m), \eqno(2)$$
where $p_{ij,k}$ are coefficients of heredity and
$$ p_{ij,k}\geq 0, \ \ \sum^m_{k=1}p_{ij,k}=1, \ \ (i,j,k=1,...,m). \eqno(3)$$
Note that each element $x\in S^{m-1}$ is a probability distribution
on $E=\{1,...,m\}.$

The population evolves by starting from an arbitrary state
(probability distribution on $E$) $x\in S^{m-1}$ then passing to the
state $Vx$ (in the next ``generation''), then to the state $V^2x$,
and so on.

For a given $x^{(0)}\in S^{m-1}$ the trajectory $\{x^{(n)}\},
n=0,1,2,...$ of $x^{(0)}\in S^{m-1}$ under action QSO (2) is defined
by $x^{(n+1)}=V(x^{(n)}),$ where $n=0,1,2,...$

 One of the main
problems in mathematical biology consists in the study of the
asymptotical behavior of the trajectories. This problem was fully
solved for Volterra QSO (see [6],[8],[12]) which is defined by (2),
(3) and the additional assumption
$$ p_{ij,k}=0, \ \ {\rm if} \ \ k\notin \{i,j\}. \eqno(4)$$
The biological treatment of condition (3) is rather clear: the
offspring repeats the genotype of one of its parents.

In paper [6] the general form of Volterra QSO $V: x=(x_1,...,x_m)\in
S^{m-1}\to V(x)=x'=(x'_1,...,x'_m)\in S^{m-1}$ is given:

$$x'_k=x_k\bigg(1+\sum_{i=1}^na_{ki}x_i\bigg), \eqno(5)$$
where $a_{ki}=2p_{ik,k}-1$ for $i\ne k$ and $a_{kk}=0.$ Moreover
$a_{ki}=-a_{ik}$ and $|a_{ki}|\leq 1.$

In papers  [6], [8] the theory of QSO (5) was  developed using
theory of the Lyapunov functions and tournaments. But non-Volterra
QSOs (i.e. which do not satisfy the condition (4)) were not in
completely studied. Because there is no any general theory which can
be applied for investigation of non-Volterra operators. To the best
of our knowledge, there are few papers devoted to such operators.
Now we shall briefly describe the history of non-Volterra operators.

In [5] it was considered the following  family of QSO $V_\l:S^2\to
S^2,$ $V_\l=(1-\l)V_0+\l V_1, 0\leq \l \leq 1,$ where
$V_0(x)=(x_1^2+2x_1x_2, x_2^2+2x_2x_3, x_3^2+2x_1x_3)$ is Volterra
operator, $V_1(x)=(x_1^2+2x_2x_3, x_2^2+2x_1x_3, x_3^2+2x_1x_2)$ is
non-Volterra operator. Note that behavior of the trajectories of
$V_0$ is very irregular ([16], [17]).

Also well known the class of bistochastic QSOs which contains
Volterra operators and non-Volterra operators as well. In [7], [13]
the class of such operators is described.

In [3] a class of kvazi-Volterra operators is introduced. For such
operators the condition (4) is not satisfied only for very few
values of $i,j,k.$

Papers [9-11] are devoted to study of non-Volterra operators which
are generated from Volterra operators (5) by a cyclic permutation of
coordinates i.e. $V_{\pi}: x'_{\pi(i)}=x_i(1+\sum^m_{k=1}a_{ik}x_k),
i=1,...,m,$ where $\pi$ is a cyclic permutation on the set of
indices $E.$

Note that each quadratic operator $V$ can be uniquely defined by a
cubic matrix ${\mathbf P}\equiv {\mathbf
P}(V)=\{p_{ij,k}\}^n_{i,j,k=1}$ with condition (3). Usually (see
e.g. [1], [3], [5]-[15], [17], [18]) the matrix ${\mathbf P}$ is
known. In [2], [4] a constructive description of ${\mathbf P}$ is
given. This construction depends on a probability measure $\mu$
which is given on a fixed graph $G$  and cardinality of a set of
cells (configurations) which can be finite or continual. In [2] it
was proven that the QSO constructed by the construction is Volterra
if and only if $G$ is a connected graph.

In  [16] using the construction of QSO for the general finite graph
and probability measure $\mu$ (here $\mu$ is product of measures
defined on maximal subgraphs of the graph $G$) a class of
non-Volterra QSOs is described. It was shown that if $\mu$ is given
by the product of the probability measures then corresponding
non-Volterra operator can be reduced to $N$ number (where $N$ is the
number of maximal connected subgraphs) of Volterra operators defined
on the maximal connected subgraphs.

Recently by the authors of this paper a new class of non-Volterra
operators is introduced. These operators satisfy
$$ p_{ij,k}=0, \ \ {\rm if} \ \ k\in \{i,j\}. \eqno(6)$$
Such an operator is called strictly non-Volterra QSO. For arbitrary
strictly non-Volterra QSO defined on $S^2$ we proved that every such
operator has unique fixed point. Also it was proven that such
operators has a cyclic trajectory and almost all trajectories
converge to the cyclic trajectory. This is quite different behavior
from the behavior of  Volterra operators, since Volterra operator
has no cyclic trajectories.

In this paper we consider an other type  of non-Volterra operators.
Which we call $F-$ quadratic stochastic operators. For a family of
such operators we show that each of them has unique fixed point and
all trajectories tend to this fixed point exponentially fast. Thus
they are ergodic operators.

\section {Definition of $F-$QSO}

In this paper we extend the set $E$ by adding element "0" i.e. we
consider $E_0=\{0,1,...,m\}.$ Fix a set $F\subset E$ and call this
set the set of "females" and the set $M=E\setminus F$ is called the
set of "males".  The element $0$ will play the role of "empty-body".

Coefficients $p_{ij,k}$ of the matrix ${\mathbf P}$ we define as
follows
$$p_{ij,k}=\left\{\begin{array}{lll}
1, \ \ {\rm if} \ \ k=0, i,j\in F\cup \{0\} \ \ {\rm or} \ \ i,j\in
M\cup \{0\};\\
0, \ \ {\rm if} \ \ k\ne 0, i,j\in F\cup \{0\} \ \ {\rm or} \ \
i,j\in
M\cup \{0\};\\
\geq 0, \ \ {\rm if } \ \ i\in F, j\in M, \forall k.
\end{array}\right.\eqno (7)$$

Biological treatment of the coefficients (7) is very clear: a
"child" $k$ can be generated if its parents are taken from different
classes $F$ and $M$. In general, $p_{ij,0}$ can be strictly positive
for $i\in F$ and $j\in M$, this corresponds, for example,  to the
case when "female" $i$ with "male" $j$ can not generate a "child"
since one of them (or both are) is ill.

\vskip 0.3 truecm

{\bf Definition 1.} {\it  For any fixed $F\subset E$, the QSO
defined by (2),(3) and (7) is called $F-$ quadratic stochastic
operator.}

\vskip 0.3 truecm

{\bf Remarks.} 1. Any $F-$ QSO is non- Volterra since $p_{ii,0}=1$
for any $i\ne 0.$

2. Note that the class of $F-$QSOs for a given $F$ does not
intersect with the classes of non-Volterra QSOs mentioned above (in
Introduction).

3. For $m=1$ there is unique $F-$QSO (independently on $F=\{1\}$ and
$F=\emptyset$) which is constant i.e. $V(x)=(1,0)$ for any $x\in
S^1.$

\section {$F-$QSO for $m=2$}

In this section we consider $m=2,$ i.e. $E_0=\{0,1,2\}.$ Take
$M=\{1\}$ and $F=\{2\}.$ We will generalize the result of this
section in the next section for $m\geq 3$ there $F$ is taken as
$F=\{2,3,...,m\}.$  Reasons of considering the case $m=2$ in a
separate section are: (i) in  case $m=2$ the number of "females" is
equal to the number of  "males"; (ii) In this case all calculations
can be done explicitly.

For $m=2,$ and $M=\{1\},$ $F=\{2\}$ the $F-$OSO is defined by matrix

$$\left(
\begin{array}{cccccc}
p_{00,0}=1 & p_{01,0}=1 & p_{02,0}=1 & p_{11,0}=1 & p_{12,0}=a & p_{22,0}=1 \\[2mm]
p_{00,1}=0 & p_{01,1}=0 & p_{02,1}=0 & p_{11,1}=0 & p_{12,1}=b & p_{22,1}=0 \\[2mm]
p_{00,2}=0 & p_{01,2}=0 & p_{02,2}=0 & p_{11,2}=0 & p_{12,2}=c & p_{22,2}=0 \\[2mm]
\end{array}
\right), $$ where
$$ a, b, c\geq 0, \ \ a+b+c=1. \eqno (8)$$
The corresponding $F-$QSO has the form

$$V_0: \left\{\begin{array}{lll}
x'_0=x_0^2+x_1^2+x_2^2+2x_0x_1+2x_0x_2+2ax_1x_2=1-2(1-a)x_1x_2,\\[3mm]
x'_1=2bx_1x_2,\\[3mm]
x'_2=2cx_1x_2.
\end{array}\right.\eqno (9)$$

The fixed point of $V$ is defined as a solution of the equation
$V(x)=x.$

The following theorem completely  describes the behavior of the
trajectories of operator (9).

\vskip 0.5 truecm

{\bf Theorem 1.} 1) {\it For any $a, b, c$ with condition (8) the
operator (9) has unique fixed point} (1,0,0).

2) {\it For any $x^{(0)}\in S^2$ the trajectory $\{x^{(n)}\}$ tends
to the fixed point $(1,0,0)$ exponentially fast.}

\vskip 0.4 truecm

{\bf Proof.} 1) It is easy to see that the solutions of the equation
$V_0(x)=x$ are only $x=e_1=(1,0,0)$ and $x=x^*=(x^*_0, x^*_1,
x^*_2)=({2bc-b-c\over 2bc}, {1\over 2c}, {1\over 2b})$ if $bc\ne 0$
and only $x=e_1$ if $bc=0.$ Note that $x^*\notin S^2.$ Indeed from
$0\leq {1\over 2b}\leq 1$ and $0\leq {1\over 2c}\leq 1$ we get
$b\geq {1\over2}$ and $c\geq {1\over 2}.$ Then by (8) we have
$b=c={1\over 2}$ i.e. $x^*_1=x^*_2=1$ this is impossible since
$x^*_0+x^*_1+x^*_2=1.$

2) For $x\in S^2$ denote $\varphi(x)=x_1x_2.$ We have
$$\varphi(x^{(n)})= \left\{\begin{array}{ll}
0, \ \ {\rm if} \ \ bc=0,\\[3mm]
(2bc)^{-1}(4bcx_1x_2)^{2^n}, \ \ {\rm if} \ \ bc\ne 0.\\[3mm]
\end{array}\right.\eqno (10)$$

Now we shall estimate $4bcx_1x_2.$ We have
$$0\leq 4bcx_1x_2\leq 4 {(b+c)^2\over 4}\cdot {(x_1+x_2)^2\over 4}\leq
{1\over 4}. \eqno (11)$$ From the second equality of (9) we get
$$x^{(n)}_1=2b\varphi(x^{(n-1)}). \eqno (12)$$
By (9)-(12) we have
$$\lim_{n\to\infty}x^{(n)}_1=\lim_{n\to\infty}x^{(n)}_2=0.$$
This completes the proof.

\section { The case $m\geq 2$}

In this section we generalize Theorem 1, but now we have not an
explicit formula for $\varphi(x^{(n)}).$

Consider $E_0=\{0,1,...,m\}, M=\{1\}, F=\{2,...,m\}.$ Then it is
easy to see that corresponding $F-$QSO has the form (cf. with (9))
$$V_1: \left\{\begin{array}{ll}
x'_0=1-2x_1\sum_{i=2}^m(1-a_{i0})x_i,\\[3mm]
x'_k=2x_1\sum^m_{i=2}a_{ik}x_i, \ \ k=1,2,...,m,\\[3mm]
\end{array}\right.\eqno (13)$$
where
$$a_{ik}=p_{1i,k}\geq 0,\ \  i=2,...,m; \ \ k=0,1,...,m; \ \ \sum_{k=0}^ma_{ik}=1,
\ \ \forall i=2,...,m. \eqno(14)$$

\vskip 0.5 truecm

{\bf Theorem 2.} {\it For any values of $a_{ik}$ with (14) the
operator (13) has unique fixed point $(1,0,0,...,0)$ (with $m$
zeros). Moreover for any $x^{(0)}\in S^m$ its trajectory
$\{x^{(n)}\}$ tends to the fixed point exponentially fast.}

\vskip 0.4 truecm

{\bf Proof.} For $x\in S^m$ denote
$$\varphi(x)=x_1\sum_{i=2}^mx_i=x_1(1-x_0-x_1). \eqno (15)$$
Using (13),(14) we get
$$x^{(n+1)}_k=2x_1^{(n)}\sum_{i=2}^ma_{ik}x_i^{(n)}\leq
2x_1^{(n)}\sum_{i=2}^mx_i^{(n)}=2\varphi(x^{(n)}), \eqno (16)$$
where $k=1,2,...,m.$

Now we shall estimate $\varphi(x^{(n+1)})$:

$$\varphi(x^{(n+1)})=x_1^{(n+1)}(1-x^{(n+1)}_0-x_1^{(n+1)})=$$
$$\bigg(2x^{(n)}_1\sum^m_{i=2}a_{i1}x^{(n)}_i\bigg)
\bigg(2x^{(n)}_1\sum^m_{i=2}(1-a_{i0})x^{(n)}_i-
2x^{(n)}_1\sum^m_{i=2}a_{i1}x^{(n)}_i\bigg). \eqno (17)$$

Using AM-GM inequality from (17) we get
$$\varphi(x^{(n+1)})\leq
4\big(x^{(n)}_1\big)^2\bigg({\sum^m_{i=2}(1-a_{i0})x^{(n)}_i\over
2}\bigg)^2\leq
\bigg(x^{(n)}_1\sum^m_{i=2}x^{(n)}_i\bigg)^2=\bigg(\varphi(x^{(n)})\bigg)^2,
\ \ n\geq 0. \eqno (18)$$ Note that $\varphi(x^{(0)})\leq {1\over
4}.$ Thus from (18) we obtain
$$\varphi(x^{(n)})\leq \Big({1\over 4}\Big)^{2^n}. \eqno(19)$$
Now it follows from (16) and (19) that
$$\lim_{n\to\infty}x^{(n)}_k=0, \ \ {\rm for\ \ any} \ \
k=1,2,...,m,$$ i.e.
$$\lim_{n\to \infty}x^{(n)}=(1,0,...,0),\ \ {\rm for\ \  any} \ \
x^{(0)}\in S^m. \eqno(20)$$

 Obviously, $(1,0,...,0)$ is a unique fixed point, since (20) holds for any
 $x^{(0)}\in S^m.$
This completes the proof.

\vskip 0.4 truecm

{\bf Remarks.} 1. The QSO $V$ satisfies the ergodic theorem if the
limit $$\lim_{n\to\infty}{1\over n}\sum^{n-1}_{j=0}x^{(j)}$$ exists
for any $x^{(0)}\in S^m.$ On the basis of numerical calculations
Ulam conjectured [15], [17], that the ergodic theorem holds for any
QSO. In [18] it was proven that this conjecture is false in general.
From Theorem 2  follows that the ergodic theorem holds for any $F-$
QSO determined by (13).

2. Assume a matrix ${\mathbf P}$ with coefficients (7) is given. Let
$N_1=N_1({\mathbf P})=|\{ij: p_{ij,0}=1\}| $ be the number of 1 in
the first row of the matrix ${\mathbf P}$, where $|A|$ denotes the
number of elements (cardinality) of $A.$  Also denote by ${\tilde
N_1}={\tilde N}_1({\mathbf P})$ the number of elements of the first
row of ${\mathbf P}$ which $<1.$ It is easy to see that ${\tilde
N}_1={|E|\over 2}(|E|+3)+1-N_1.$ By (7) one can see that
$$N_1\geq {1\over 2}\Big(|F|^2+|M|^2+3|E|\Big)+1,$$
$${\tilde N}_1\leq |F|\cdot |M|.$$ Obviously,
$${1\over 2}\Big(|F|^2+|M|^2+3|E|\Big)+1>|F|\cdot |M|.$$
Hence $N_1>{\tilde N}_1,$ for any $F\subset E$ and any ${\mathbf P}$
with elements  (7).  By this property of the matrix ${\mathbf P}$
one can say that $0$ has more priority than other elements of $E_0.$
Thus we can make a

 \vskip 0.4 truecm

{\bf Conjecture.} {\it An analogue of Theorem 2 is true for any
$F-$QSO i.e. for every $F\subset E,$ and ${\mathbf P}$ with
elements} (7).

\vskip 0.4 truecm

3. The $F-$QSOs with $F\ne \{2,3,...,m\}$ will be considered in a
separate paper.

 \newpage

 {\bf Acknowledgments.} A part of this work was done
within the scheme of Junior Associate at the ICTP, Trieste, Italy
and the first author (UAR) thanks ICTP for providing financial
support and all facilities (in May - August 2006). The final part of
this work was done at  the IHES, Bures-sur-Yvette, France. UAR
thanks the IHES for support and kind hospitality (in October -
December 2006). We also gratitude to professor R. N. Ganikhodzhaev
for many helpful discussions.

\vskip 0.4 truecm

{\bf References}

1. Bernshtein S.N.,  Solution of a mathematical problem connected
with the theory of heredity, {\it Uch. Zap. Nauchno-Issled. kaf.
Ukr. Otd. Mat.,} {\bf 1} : 83-115 (1924).

2. Ganikhodjaev N.N., An application of the theory of Gibbs
distributions to mathematical genetics, {\it Doklady Math.} {\bf
61}: 321-323 (2000).

3. Ganikhodjaev N.N., Mukhitdinov R.T.,  On a class of non-Volterra
quadratic operators, {\it Uzbek Math. Jour.} No. 3-4: 65-69 (2003).

4. Ganikhodjaev N.N., Rozikov U.A.,  On quadratic stochastic
operators generated by Gibbs distributions, {\it Regular and Chaotic
Dynamics.} {\bf 11}: No. 3 (2006).

5. Ganikhodzhaev, R. N. A family of quadratic stochastic operators
that act in $S^2$. {\it Dokl. Akad. Nauk UzSSR.}  No. 1: 3-5 (1989).

6. Ganikhodzhaev R.N., Quadratic stochastic operators, Lyapunov
functions and tournaments, {\it Russian Acad. Sci. Sbornik Math.}
{\bf 76}: 489-506 (1993).

7. Ganikhodzhaev, R. N. On the definition of quadratic bistochastic
operators. {\it Russian Math. Surveys.} {\bf 48}, no. 4: 244-246
(1992).

8. Ganikhodzhaev R.N., A chart of fixed points and Lyapunov
functions for a class of discrete dynamical systems. {\it Math.
Notes} {\bf 56}: 1125-1131 (1994).

9. Ganikhodzhaev, R. N.; Dzhurabaev, A. M. The set of equilibrium
states of quadratic stochastic operators of type $V_{\pi}$. {\it
Uzbek. Mat. Zh.}  No. 3: 23-27 (1998).

10. Ganikhodzhaev, R. N.; Karimov, A. Z. On the number of vertices
of a polyhedron of bistochastic quadratic operators. {\it Uzbek.
Mat. Zh.} No. 6: 29-35 (1999).

11. Ganikhodzhaev, R. N.; Abdirakhmanova, R. E. Description of
quadratic automorphisms of a finite-dimensional simplex. {\it Uzbek.
Mat. Zh.}  No. 1: 7--16 (2002).

12. Ganikhodzhaev, R.N., Eshmamatova D.B. Quadratic automorphisms of
simplex and asymptotical behavior of their trajectories. {\it
Vladikavkaz Math. Jour.} {\bf 8}: 12-28 (2006).

13. Ganikhodzhaev, R. N.; Eshniyazov, A. I. Bistochastic quadratic
operators. {\it Uzbek. Mat. Zh.} No. 3: 29-34 (2004).

14. Kesten H. Quadratic transformations: a model for population
growth I, II. {\it Adv. Appl. Prob.} No.2: 1-82 and 179-228 (1970).

15. Lyubich Yu.I. Mathematical structures in population genetics.
{\sl Biomathematics}, {\bf 22,} Springer-Verlag, 1992.

16. Rozikov U.A., Shamsiddinov N.B. On non-Volterra quadratic
stochastic operators generated by a product measure. {\it ICTP
preprint 2006}, and arXiv:math.DS/0608201.

17. Stein P.R., Ulam S.M. Nonlinear transformations studies on
electronic computers. {\it Rozprawy Mat.} {\bf 39}: 1-15 (1964).

18. Zakharevich M.I., The behavior of trajectories and the ergodic
hypothesis for quadratic mappings of a simplex. {\it Russian Math.
Surveys}, {\bf 33}: 207-208 (1978).

\end{document}